\documentclass[11pt]{article}

\usepackage{amsmath}    % need for subequations
\usepackage{graphicx}   % need for figures
\usepackage{verbatim}   % useful for program listings
\usepackage{color}      % use if color is used in text
\usepackage{subfigure}  % use for side-by-side figures
\usepackage{amssymb,tikz}
\usepackage{soul}

\newtheorem{theorem}{Theorem}[section]
\newtheorem{proposition}[theorem]{Proposition}

\newtheorem{conjecture}[theorem]{Conjecture}
\newtheorem{corollary}[theorem]{Corollary}
\newtheorem{lemma}[theorem]{Lemma}

\newtheorem{problem}[theorem]{Problem}

\newcommand{\proof}{\noindent{\bf Proof.\ }}
\newcommand{\qed}{\hfill $\square$\medskip}
\def\cp{\,\square\,}

\newcommand{\N}{\mathbb N}

\DeclareMathOperator {\sg} {sg}

\DeclareMathOperator {\diam} {diam}

\let\d\relax
\DeclareMathOperator {\d} {d}
\DeclareMathOperator {\n} {n}

\begin{document}

% Define style for nodes
\tikzstyle{every node}=[circle, draw, fill=black!10,
                        inner sep=0pt, minimum width=4pt]

\title{Strong geodetic problem on Cartesian products of graphs}
 
\author{
	Vesna Ir\v si\v c $^{a}$
	\and
	Sandi Klav\v zar $^{a,b,c}$
}

\date{August 5, 2017}

\maketitle
% \vspace{-0.8 cm}
\begin{center}
	$^a$ Institute of Mathematics, Physics and Mechanics, Ljubljana, Slovenia\\
	\medskip

	$^b$ Faculty of Mathematics and Physics, University of Ljubljana, Slovenia\\
% 	{\tt sandi.klavzar@fmf.uni-lj.si}\\
	\medskip
	
	$^c$ Faculty of Natural Sciences and Mathematics, University of Maribor, Slovenia
		
\end{center}

\begin{abstract}
The strong geodetic problem is a recent variation of the geodetic problem. For a graph $G$, its strong geodetic number $\sg(G)$ is the cardinality of a smallest vertex subset $S$, such that each vertex of $G$ lies on a fixed shortest path between a pair of vertices from $S$. In this paper, the strong geodetic problem is studied on the Cartesian product of graphs. A general upper bound for $\sg(G \cp H)$ is determined, as well as exact values for $K_m \cp K_n$, $K_{1, k} \cp P_l$, and certain prisms. Connections between the strong geodetic number of a graph and its subgraphs are also discussed.
\end{abstract}


\medskip
\noindent{\bf AMS Subj.\ Class.: 05C12, 05C70, 05C76; 68Q17}

%%%%%%%%%%%%%%%%%%%%%%%%%%%
\section{Introduction}
%%%%%%%%%%%%%%%%%%%%%%%%%%%

Covering vertices of a graph with shortest paths is a natural (optimization) problem arising from different applied problems that respectively led to several different graph theory models. The seminal of them, the {\em geodetic problem}~\cite{HLT93}, aims to find a smallest subset of vertices of a given graph such that the geodesics between them cover all its vertices, see the review~\cite{BrKo11}. Recent studies on this problem have focused on characterizations of graphs with large geodetic number~\cite{Filomat}, on geodesic graphs~\cite{soloff-2015}, and on connections between the geodetic problem and a block decomposition~\cite{ekim-2014}. Applications of the geodetic problem can be found in convexity theory~\cite{Centeno, Jiang+2004, Lu2004, pelayo-2013} and in  game theory~\cite{igre}.

Another variation of the problem of covering vertices with shortest paths is the {\em isometric path problem}~\cite{fisher-2001} where the aim is to determine the minimum number of shortest paths required to cover all the vertices of a graph. Following~\cite{fisher-2001} this problem has been investigated on Cartesian products of graphs~\cite{Fitzpatrick-1999-CN}, in particular on Hamming graphs as well as on complete $r$-partite graphs in~\cite{pan-2006}.

Motivated by applications in social networks, the strong geodetic problem was introduced in~\cite{MaKl16a} as follows. Let $G=(V,E)$ be a graph. Given a set $S\subseteq V$, for each pair of vertices $\{x,y\}\subseteq S$, $x\ne y$, let $\widetilde{g}(x,y)$ be a {\em selected fixed} shortest path between $x$ and $y$. We set  
$$\widetilde{I}(S)=\{\widetilde{g}(x, y) : x, y\in S\}\,,$$ and $V(\widetilde{I}(S))=\bigcup_{\widetilde{P} \in \widetilde{I}(S)} V(\widetilde{P})$. If $V(\widetilde{I}(S)) = V$ for some $\widetilde{I}(S)$, then the set $S$ is called a {\em strong geodetic set}. For a graph $G$ with just one vertex, we consider the vertex as its unique strong geodetic set. The {\em strong geodetic problem} is to find a minimum strong geodetic set of $G$. The cardinality of a minimum strong geodetic set is the {\em strong geodetic number} of $G$ and is denoted by $\sg(G)$. 

In the first paper~\cite{MaKl16a} on the strong geodetic number this invariant  has been determined for complete Apollonian networks and proved that the problem is NP-complete. Then, in~\cite{Klavzar+2017}, the problem was studied on grids and cylinders. Among other results it was proved that if $r$ is large enough comparing to $n$, then $\sg(P_r\cp P_n)= \lceil 2\sqrt{n}\, \rceil$. Some general properties of the strong geodesic problem, in particular with respect to the diameter, and a solution for balanced complete bipartite graphs has been very recently reported in~\cite{bipartite}. We also refer to~\cite{MaKl16b} for an edge version of the problem.

In this paper, the strong geodesic problem is studied on Cartesian product graphs. In the next section we give several upper bounds on $\sg(G \cp H)$ and study their sharpness. In Section~\ref{sec:examples} we determine the strong geodetic number for several families of Cartesian products, including products of complete graphs. We also discuss a possible lower bound for $\sg(G \cp H)$. Motivated by this discussion, in the final section we focus on possible connections between the strong geodetic number of a graph and its subgraphs. But first we list necessary definitions. 

All graphs considered in this paper are simple and connected. The {\em distance} $d_G(u,v)$ between vertices $u$ and $v$ of a graph $G$ is the number of edges on a shortest $u,v$-path ($u,v$-{\em geodesic}). The {\em diameter} ${\rm diam}(G)$ of $G$ is the maximum distance between vertices of $G$. We denote the order of a graph by $\n(G)$. A vertex $v$ of a graph $G$ is {\em simplicial} if its neighborhood induces a clique. We will use the notation $[n] = \{1,\ldots, n\}$ and the convention that $V(P_n) = V(K_n) = V(C_n) = [n]$ for any $n \geq 1$, where the edges of the path $P_n$, the complete graph $K_n$, and the cycle $C_n$ are defined in the natural way.

The {\em Cartesian product} $G\cp H$ of graphs $G$ and $H$ is the graph with vertices $V(G) \times V(H)$, where the edges $(g,h)$ and $(g',h')$ are adjacent if either $g=g'$ and $hh'\in E(H)$, or $h=h'$ and $gg'\in E(G)$. If $h\in V(H)$, then a subgraph of $G\cp H$ induced by the set of vertices $\{(x,h) \; ; \; x\in V(G)\}$ is isomorphic to $G$;  it is denoted by $G^h$ and called a {\em $G$-layer}, a {\em horizontal layer} or a {\em row}. Analogously $H$-layers are defined; if $g\in V(G)$, then the corresponding \emph{$H$-layer}, called a {\em vertical layer} or a {\em column}, is denoted $^{g}H$.

%%%%%%%%%%%%%%%%%%%%%%%%%%%%%%%%%%%%%%%%%%%%%%%%%%%%%%%%%%
%%%%%%%%%%%%%%%%%%%%%%%%%%%%%%%%%%%%%%%%%%%%%%%%%%%%%%%%%%
\section{Upper bounds on $\sg(G \cp H)$}
\label{sec:generalUpper}
%%%%%%%%%%%%%%%%%%%%%%%%%%%%%%%%%%%%%%%%%%%%%%%%%%%%%%%%%%
%%%%%%%%%%%%%%%%%%%%%%%%%%%%%%%%%%%%%%%%%%%%%%%%%%%%%%%%%%

The investigations from~\cite{Klavzar+2017} indicate that it is not easy to determine the strong geodetic number of an arbitrary integer grid, that is, $\sg(P_r\cp P_n)$. As these grids are among the simplest Cartesian product graphs, it would be too ambitious to expect a formula for $\sg(G\cp H)$. In this section we therefore consider upper bounds for $\sg(G\cp H)$ and discuss their sharpness.

Note first that lifting a strong geodetic set of $G$ (resp.\ $H$) into each of the $G$-layers (resp.\ $H$-layers) yields $\sg(G \cp H) \leq \min\{ \sg(G) \n(H), \sg(H) \n(G)\}$. This observation can be improved as follows.

\begin{theorem}
\label{thm:generalUpper}
If $G$ and $H$ are graphs, then
$$\sg(G \cp H) \leq \min\{ \sg(H) \n(G) - \sg(G) + 1, \sg(G) \n(H) - \sg(H) + 1 \}.$$
\end{theorem}

\proof
Since the Cartesian product operation is commutative, it suffices to prove that  $\sg(G \cp H) \leq \sg(H) \n(G) - \sg(G) + 1$. 

Let $S_G$ be a strong geodetic set of $G$, $\widetilde{I}(S_G)$ fixed geodesics in $G$, $S_H$ a strong geodetic set of $H$, and $\widetilde{I}(S_H)$ fixed geodesics in $H$, where $|S_G| = \sg(G) = k$ and $|S_H| = \sg(H) = l$. Set $S_G = \{ g^0, g^1, \ldots, g^{k-1} \}$ and $S_H = \{ h^0, h^1, \ldots, h^{l-1} \}$. Denote with $P_i$ the $g^0, g^i$-geodesic from $\widetilde{I}(S_G)$ for all $i \in [k-1]$ and with $Q_j$ the $h^0, h^j$-geodesic from $\widetilde{I}(S_H)$ for all $j \in [l-1]$. 

Define $T = (V(G) \times S_H) - \{ (g, h^0) \; ; \; g \in S_G - \{ g^0 \} \}$. Clearly, $|T| = \sg(H) \n(G) - \sg(G) + 1$. We claim that $T$ is a strong geodetic set of $G\cp H$. To show it, we first fix geodesics in $H$-layers between vertices from $T$ in the same way as they are fixed in $\widetilde{I}(S_H)$. The only (possibly) uncovered vertices are the ones lying in $H$-layers $^{g^i}H$ for $i \in [k-1]$ that lie on paths $Q_j$ for $i \in [l-1]$ and $j \in [l-1]$. To cover them we fix $(g^i,h^j), (g^0, h^0)$-geodesics as paths $\{g^i\} \times Q_j$ joined with $P_i \times \{ h^0 \}$ for all $i \in [k-1], j \in [l-1]$. In this way all the vertices of $G \cp H$ are covered, hence $\sg(G \cp H) \leq |T|$.  
\qed

If $n\ge 2$, then $\sg(P_n \cp K_2) = 3 = \sg(P_n) \n(K_2) - \sg(K_2) + 1$. This example shows that the inequality of Theorem~\ref{thm:generalUpper} is best possible. To construct more sharpness examples we need the following general property.

\begin{lemma}
\label{lem:simplicial}
If $G$ and $H$ are graphs, $v$ is a simplicial vertex of $G$, and $S$ is a strong geodetic set of $G \cp H$, then $S\cap\, ^{v}H \ne \emptyset$. 
\end{lemma}

\proof
Suppose on the contrary that $S\cap\, ^{v}H = \emptyset$. Let $P\in \widetilde{I}(S)$ be an arbitrary geodesic that contains some vertices of $^{v}H$. By the assumption, $P$ starts and ends outside $^{v}H$. Let $(g,h)$ be the first vertex of $P$ with a neighbor in $^{v}H$ and let $(g',h')$ be the first subsequent vertex of $P$ that does not lie in $^{v}H$. Suppose $g\ne g'$. Then a $((g,h),(g,h'))$-geodesic together with the edge $(g,h')(g',h')$ (which exists since $v$ is a simplicial vertex of $G$) yields a shorter $((g,h),(g',h'))$-path than the $((g,h),(g',h'))$-subpath of $P$, a contradiction with the fact that $P$ is a geodesic. If $g=g'$ we get the same contradiction, except that there is no need to add the edge $(g,h')(g',h')$.  
\qed

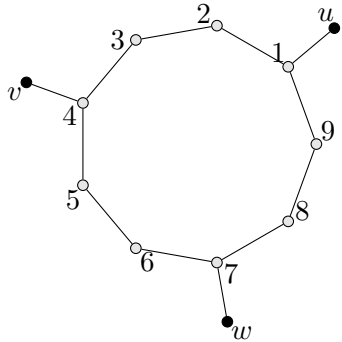
\begin{figure}[!ht]
\begin{center}
    \begin{tikzpicture}[scale=0.8]
    \pgfmathtruncatemacro{\N}{9}
    \pgfmathtruncatemacro{\r}{2}
    \begin{scope}
        \foreach \x in {1,...,\N}
            \node[label={60 + \x*360/\N:$\x$}] (\x) at (\x*360/\N:\r cm) {};
        
        \foreach \x [remember=\x as \lastx (initially 1)] in {1,...,\N,1}
            \path (\x) edge (\lastx);
            
        \node[label={60 + 1*360/\N:$u$}, fill=black] (u) at (1*360/\N:3 cm) {};
        \path (1) edge (u);
        \node[label={60 + 4*360/\N:$v$}, fill=black] (u) at (4*360/\N:3 cm) {};
        \path (4) edge (u);
        \node[label={60 + 7*360/\N:$w$}, fill=black] (u) at (7*360/\N:3 cm) {};
        \path (7) edge (u);
    \end{scope}
    \end{tikzpicture}
    \caption{A graph $G_3$ and its strong geodetic set.}
    \label{fig:cycle}
\end{center}
\end{figure}

If $n\ge 3$, then let $G_n$ be the graph obtained from $C_{3n}$ by adding vertices $u, v, w$ and edges $u \sim 1$, $v \sim n+1$, and $w \sim 2n+1$; cf.~Fig.~\ref{fig:cycle} for $G_3$. 

% \iffalse
% \fi

Recall from~\cite{MaKl16a} that a simplicial vertex lies in every strong geodetic set. Hence $\sg(G_n) \geq 3$. On the other hand, $\{u, v, w\}$ is a strong geodetic set which implies that $\sg(G_n) = 3$. Consider now the product $G_n \cp K_2$. By Lemma~\ref{lem:simplicial} we have $\sg(G_n \cp K_2) \geq 3$. Suppose $\sg(G_n \cp K_2) = 4$ and let $S$ be a strong geodetic set with $|S|=4$. Then $S$ must have two vertices in each of the $G_n$-layers. Thus, applying Lemma~\ref{lem:simplicial} again, we can assume without loss of generality that $\{(u, 1), (v, 2), (w, 2)\} \subseteq S$. If $s$ is the fourth vertex of $S$, then $s$ lies in $G_n^1$ and equals one of $(n+2, 1), \ldots, (2n, 1)$, for otherwise these vertices could not lie on any geodesic from $\widetilde{I}(S)$. Without loss of generality assume that the $s, (u, 1)$-geodesic passes the vertex $(n+1, 1)$. But then it is not possible to cover all vertices $(2n+2, 1), (2n+2, 2), \ldots, (2n-1, 1), (2n-1, 2)$, as $n \geq 3$. In conclusion, 
$$\sg(G_n \cp K_2) = 5 = \sg(G_n) \n(K_2) - \sg(K_2) + 1\,,$$
hence we have constructed another infinite family attaining equality in Theorem~\ref{thm:generalUpper}.

If $G=(V,E)$ is a graph and $S\subseteq V$, then $S$ is called a {\em $2$-packing} if $d(x,y)\ge 3$ holds for any $x,y\in S$, $x\ne y$. Equivalently, $S$ is not a $2$-packing if and only if $S$ contains vertices $u\ne v$ such that $\d(u,v) \leq 2$. Now we can improve Theorem~\ref{thm:generalUpper} in the following case.  

\begin{proposition}
\label{prop:-2}
If $G$ is a graph with $\sg(G) \geq 3$ that admits a strong geodetic set $S$ which is not a $2$-packing, then 
$$\sg(G \cp K_2) \leq 2 \sg(G) - 2.$$
\end{proposition}

\proof
Let $S$ be a strong geodetic set of a graph $G$ with the desired properties: $|S| = \sg(G) = k$ and $S = \{ u, v, u_1, \ldots, u_{k-2} \}$, where $d(u,v)\le 2$. Let $\widetilde{I}(S)$ be a set of fixed geodesics. Let $P_{uv} \in \widetilde{I}(S)$ be the path between $u$ and $v$ and note that the length of $P_{uv}$ is either $1$ or $2$. For $i \in [k-2]$ denote by $P_i \in \widetilde{I}(S)$ the $u, u_i$-geodesic and by $Q_i \in \widetilde{I}(S)$ the $v, u_i$-geodesic. 

Set $T = ((S - \{ u\}) \times \{1\}) \cup ((S - \{v\}) \times \{2\})$. Clearly, $|T| = 2 \sg(G) - 2$. Fix the same geodesics as in $\widetilde{I}(S)$ between vertices in $(S - \{ u \}) \times \{1\}$ and between vertices in $(S - \{ v \}) \times \{2\}$. The only possibly uncovered vertices are the ones lying on paths $P_i$ in $G^1$ and on paths $Q_i$ in $G^2$. Thus we also fix geodesics $P_i \times \{1\}$ in $G^1$ joined with an edge $(u, 1) \sim (u, 2)$ for all $i \in [k-2]$ and geodesics $Q_i \times \{2\}$ in $G^2$ joined with $(v, 2) \sim (v, 1)$ for all $i \in [k-2]$. 

If $\d(u, v) = 1$, all vertices are already covered. If $\d(u, v) = 2$ and $w$ is the remaining vertex on the path $P_{uv}$, the only uncovered vertices are $(w, 1)$ and $(w, 2)$. These two remaining uncovered vertices can be covered with the geodesic $(v,1) \sim (w,1) \sim (w,2) \sim (u,2)$. Hence, $\sg(G \cp K_2) \leq 2 \sg(G) - 2$.
\qed

\begin{corollary}
\label{cor:-2}
If $G$ is a graph with $\diam(G) = 2$ and $\sg(G) \geq 3$, then 
$$\sg(G \cp K_2) \leq 2 \sg(G) - 2.$$
\end{corollary}

We point out that Proposition~\ref{prop:-2} and Corollary~\ref{cor:-2} do not hold in the case when $\sg(G) = 2$. For instance, if $G=K_2$, then $\sg(K_2 \cp K_2) = 3 \nleq 2 = 2 \sg(K_2) - 2$. 

Using a reasoning parallel to the one from the proof of Proposition~\ref{prop:-2}, the following generalization can be derived.

\begin{proposition}
\label{prop:-n}
If $G$ is a graph with $\sg(G) \geq 3$ that admits a strong geodetic set which is not a $2$-packing, then
$$\sg(G \cp K_n) \leq n \sg(G) - n.$$
\end{proposition}

Based on the above ideas, we can state our second main result of this section that generalizes Proposition~\ref{prop:-2} and in a special case decreases by $1$ the upper bound of Theorem~\ref{thm:generalUpper}. 

\begin{theorem}
\label{thm:-2general}
If $G$ is a graph, and $H$ is a graph with $\sg(H) \geq 3$ that admits a strong geodetic set which is not a $2$-packing, then 
$$\sg(G \cp H) \leq \sg(H) \n(G) - \sg(G).$$
\end{theorem}

\proof
Let $S_H$ be a strong geodetic set of a graph $H$ with the desired properties: $|S| = \sg(H) = l \ge 3$ and $S_H = \{ u, v, h^1, \ldots, h^{l-2}\}$, where $d(u,v)\le 2$. Let $\widetilde{I}(S_H)$ be a set of fixed geodesics that cover $V(H)$. Let $P_{uv} \in \widetilde{I}(S_H)$ be the path between $u$ and $v$. Denote by $P_i \in \widetilde{I}(S_H)$ a fixed $u, h^i$-geodesics and by $Q_i \in \widetilde{I}(S_H)$ a fixed $v, h^i$-geodesics for all $i \in [l-2]$. 

Let $S_G$ be a strong geodetic set of $G$, $\widetilde{I}(S_G)$ fixed geodesics and $|S_G| = \sg(G) = k$. Set $S_G = \{ w, g^1, \ldots, g^{k-1} \}$. Denote with $R_i$ a fixed $w, g^i$-geodesic from $\widetilde{I}(S_G)$ for all $i \in [k-1]$.

Set $T = (V(G) \times S_H) - \{ (g, u) \; ; \; g \in S_G - \{ w \} \} - \{ (w, v) \}$. Clearly, $|T| = \sg(H) \n(G) - \sg(G)$. Geodesics in $H$-layers between vertices from $T$ are fixed in the same way as in $\widetilde{I}(S_H)$. The only (possibly) uncovered vertices are the ones lying in $H$-layers $^{g^i}H$ for $i \in [k-1]$ that lie on paths $P_j$ for $j \in [l-2]$ and those on paths $Q_j$ in the layer $^{w}H$ for $j \in [l-2]$. Thus we also fix $(g^i,h^j), (w, u)$-geodesics as paths $\{g^i\} \times P_j$ joined with $R_i \times \{ u \}$ for all $i \in [k-1], j \in [l-2]$ and $(w, h^j), (g^i, v)$-geodesics as paths $\{w\} \times Q_j$ joined with $R_i \times \{v\}$ for all $i \in [k-1], j \in [l-2]$. 

If $\d(u, v) = 1$, all vertices of $G \cp H$ are already covered. If $\d(u, v) = 2$ and $t$ is the remaining vertex on $P_{uv}$, then we also fix geodesic $(g^i, v) \sim (g^i, t) \sim (w, t) \sim (w, u)$ for all $i \in [k-1]$. Now all vertices of $G \cp H$ are covered, hence $\sg(G \cp H) \leq |T|$.
\qed

\begin{corollary}
\label{cor:-2general}
If $G$ and $H$ are graphs with $\diam(G) = \diam(H) = 2$ and $\sg(G), \sg(H) \geq 3$, then 
$$\sg(G \cp H) \leq \min\{  \sg(H) \n(G) - \sg(G), \sg(G) \n(H) - \sg(H) \}.$$
\end{corollary}

%%%%%%%%%%%%%%%%%%%%%%%%%%%%%%%%%%%%%%%%%%%%%%%%%%%%%%%%%%
%%%%%%%%%%%%%%%%%%%%%%%%%%%%%%%%%%%%%%%%%%%%%%%%%%%%%%%%%%
\section{Exact values for some Cartesian products}
\label{sec:examples}
%%%%%%%%%%%%%%%%%%%%%%%%%%%%%%%%%%%%%%%%%%%%%%%%%%%%%%%%%%
%%%%%%%%%%%%%%%%%%%%%%%%%%%%%%%%%%%%%%%%%%%%%%%%%%%%%%%%%%

In this section we determine the strong geodetic number of certain prisms (Theorem~\ref{thm:equality_prisms}), of $K_{1,k} \cp P_l$ (Proposition~\ref{prop:starPn}), and of Hamming graphs $K_m \cp K_n$ (Theorem~\ref{thm:KmKn}). At the end of the section we pose a conjecture asserting a general lower bound on $\sg(G\cp H)$. The conjecture has been verified for small prisms by computer and is, provided it holds true, best possible by the results of this section.

\begin{theorem}
\label{thm:equality_prisms}
(i) If $n \geq 5$ is an integer, then $\sg(K_n - e) = \sg((K_n - e) \cp K_2) = n-1$.

(ii) If $G$ is a graph, $S$ the set of its simplicial vertices, $|S| \geq 4$, and $S$ is a strong geodetic set of $G$, then $\sg(G \cp K_2) = \sg(G)$.
\end{theorem}

\proof
(i) Let $G = K_n - e$ and $e = \{u, v\}$, $u \nsim v$. Denote $V(G) = \{u, v, x_1, \ldots, x_{n-2}\}$. As $G$ is not a complete graph, it follows from~\cite{bipartite} that $\sg(G) \leq n-1$. Let $S$ be a minimum strong geodetic set of $G$. As vertices $u$ and $v$ are simplicial, $u, v \in S$. Any $u,v$-geodesic covers exactly one other vertex, say $x_{n-2}$. Thus $S - \{u, v\}$ is a strong geodetic set of $G - \{u, v, x_{n-2}\}$, a complete graph on $n-3$ vertices. Hence, $\sg(G) \geq 2 + \sg(K_{n-3}) = n-1$. 

We now prove that $\sg(G \cp K_2) \leq n-1$. Consider $S = \{ (u, 1), (u, 2), (v, 1), (v,2) \}$ and $T = \{(x_i, 1) \; ; \; i \in \{ 4,\ldots, n-2 \}\}$. Geodesics between vertices from $S$ can be fixed in such a way, that $\{ (x_i, j) \; ; \; i \in [3], j \in [2] \}$ are all covered. The remaining uncovered vertices can be covered with geodesics $(x_i, 1) \sim (x_i, 2) \sim (u, 2)$. Hence, $S \cup T$ is a strong geodetic set of a graph $G \cp K_2$ and $\sg(G \cp K_2) \leq |S \cup T| = 4 + (n-5) = n-1$.

It remains to prove that $\sg(G \cp K_2) \geq n-1$. Notice that the longest geodesics and the only ones of length $3$ in graph $G \cp K_2$ are $(u,1), (v,2)$- and $(u,2), (v,1)$-geodesics. All other geodesics are of length $1$ or $2$ and can therefore cover at most one $K_2$-layer. Furthermore, any $K_2$-layer that is not covered with one of the longest geodesics must contain at least one vertex from the strong geodetic set. Let $S$ be the minimum strong geodetic set of $G \cp K_2$ and $\widetilde{I}(S)$ the fixed geodesics. Consider the following cases.
\begin{enumerate}
\item[(a)] If $\widetilde{I}(S)$ contains two longest geodesics, then geodesics between vertices $\{ (u, 1),$ $(u, 2),$ $(v, 1),$ $(v,2) \}$ can cover five different $K_2$-layers. To cover the remaining $n-5$ $K_2$-layers, $S$ must contain at least $n-5$ more vertices. Hence, $|S| \geq n-1$.
\item[(b)] If $\widetilde{I}(S)$ contains only one of the longest geodesics, this geodesic lies in  three $K_2$-layers. To cover the remaining $n-3$ $K_2$-layers, we need at least $n-3$ more vertices. Hence, $|S| \geq 2 + (n-3) = n-1$.
\item[(c)] If $\widetilde{I}(S)$ contains none of the longest geodesics, then at most one vertex among $\{ (u, 1), (u, 2), (v, 1), (v,2) \}$ lies in $S$. Thus at least $n-1$ $K_2$-layers are still completely uncovered, hence $|S| \geq n-1$. 
\end{enumerate}
It follows from the above, that $\sg(G \cp K_2) \geq n-1$.

\medskip
(ii) Clearly, $\sg(G) = |S|$. By Lemma~\ref{lem:simplicial} we have $\sg(G \cp K_2) \geq |S|$. Now we prove that the equality is attained. 

Let $k = |S|$, $\widetilde{I}(S)$ fixed geodesics that cover $G$, and $S = \{l_1, \ldots, l_{k-2}, r_1, r_2\}$. For $s, t \in S$ denote the fixed $s, t$-geodesic from $\widetilde{I}(S)$ by $s \leadsto t$. Set $T = \{ (l_1, 2), \ldots, (l_{k-2}, 2),$ $(r_1, 1), (r_2, 1) \}$. Fix geodesics between vertices from $T$ as follows:
$$(l_i, 2) \leadsto (r_{e(i)}, 2) \sim (r_{e(i)}, 1),$$
where $$e(i) = \begin{cases}
i; & i \in [2],\\
2; & \text{otherwise},
\end{cases}$$
and 
$$(l_i, 2) \sim (l_i, 1) \leadsto (r_{f(i)}, 1),$$
where
$$f(i) = \begin{cases}
2; & i =1,\\
1; & \text{otherwise}.
\end{cases}$$

Geodesics of the first type cover all vertices of the form $(u, 2)$, where $u \in V(G)$, and the second type covers all vertices $(u, 1)$, where $u \in V(G)$. Hence, $T$ is a strong geodetic set and $\sg(G \cp K_2) = |T| = |S| = \sg(G)$. 
\qed

With respect to Theorem~\ref{thm:equality_prisms}(ii) we note that for the geodetic problem the graphs $G$ with the property that the set of simplicial vertices of $G$ is geodetic, were studied under the name {\em extreme geodesic graphs}~\cite{chartrand-2002}. Notice also that Theorem~\ref{thm:equality_prisms}(i) does not hold for $n \leq 4$, as $\sg(K_4 - e) = 3 < 4 = \sg((K_4 - e) \cp K_2)$. Moreover, the products $P_n \cp K_2$ and $K_{3,1} \cp K_2$ demonstrate that Theorem~\ref{thm:equality_prisms}(ii) does not hold for $|S| \leq 3$. 

We now derive two exact results for Cartesian products which are not prisms. The first one reads as follows (and is in a way a generalisation of Theorem~\ref{thm:equality_prisms}(ii)).

\begin{proposition}
\label{prop:starPn}
If $k, l$ are integers, $k \geq 5$ and $l \geq 1$, then $\sg(K_{1, k} \cp P_l) = \sg(K_{1, k})$.
\end{proposition}

\proof
The graph $K_{1,k}$ is a tree with $k$ leaves, hence $\sg(K_{1, k}) = k$ and $\sg(K_{1,k} \cp P_l) \geq k$. 

Let $V(K_{1, k}) = \{ v, l_1, \ldots, l_{k-2}, r_1, r_2 \}$ where $v$ is the vertex of degree $k$. Define $S = \{ (l_1, l), \ldots, (l_{k-2}, l) , (r_1, 1), (r_2, 1) \}$. As shortest paths in $P_l$ are unique, $x, y$-geodesic can be denoted by $x \leadsto y$. Fix geodesics between vertices from $S$ in the following way:
$$(l_i, l) \sim (v, l) \sim (r_i, l) \leadsto (r_i, 1)$$
 for $i \in [2]$, 
$$(l_i, l) \sim (v, l) \leadsto (v, 1) \sim (r_2, 1)$$
 for $i \in \{3, \ldots, k-2\}$, 
and 
$$(l_i, l) \leadsto (l_i, 1) \sim (v, 1) \sim (r_{f(i)}, 1),$$
where
$$f(i) = \begin{cases}
2; & i = 1,\\
1; & i \ne 1.
\end{cases}$$

Clearly, these geodesics cover all vertices of the graph (as $k-2 \geq 3$), hence $\sg(K_{1, k} \cp P_l) = k$.
\qed

Proposition~\ref{prop:starPn} does not hold for $k \leq 4$ if $l \geq 3$ (the cases $l \in \{1, 2\}$ are simple). Consider the following example. Let $V(K_{1, 4}) = \{v, l_1, l_2, r_1, r_2\}$ as above. Suppose $\sg(K_{1, 4} \cp P_l) = 4$. If $T^l$ (or equivalently $T^1$) contains only one vertex from a minimum strong geodetic set, say $l_1$, then geodesics from $(l_1, l)$ to the other three vertices must pass vertices $(l_2, l), (r_1, l), (r_2, l), (l_1, l-1)$ which is not possible. Hence, any strong geodetic set of size $4$ contains two vertices in the layer $T^1$ and two vertices in $T^l$. Without loss of generality let $S = \{ (l_1, l), (l_2, l), (r_1, 1), (r_2, 1) \}$ be a minimum strong geodetic set. Geodesics $(l_1, l) \sim (v, l) \sim (l_2, l)$ and $(r_1, 1) \sim (v, 1) \sim (r_2, 1)$ are clearly fixed. Each of the remaining four geodesics can cover at most $l-1$ uncovered vertices. But the graph has $4(l-1) + (l-2)$ vertices to cover, hence $\sg(K_{1,4} \cp P_l) \geq 5$. Since the set $S \cup \{  (v, 1) \}$ is a strong geodetic set, we have $\sg(K_{1,4} \cp P_l) = 5$.  

Our last exact result is the following.

\begin{theorem}
\label{thm:KmKn}
If $m, n$ are positive integers and $m \geq n$, then
$$\sg(K_m \cp K_n) = \begin{cases}
2 n - 1; & m = n,\\
2n; & n < m < 2n,\\
m; & m  \geq 2n.
\end{cases}$$
\end{theorem}

\proof
Since every vertex of a complete graph is simplicial, Lemma~\ref{lem:simplicial} implies that any strong geodetic set of $K_m \cp K_n$ contains at least one vertex from each row and at least one vertex from each column, hence $\sg(K_m \cp K_n) \geq \max\{m, n\} = m$. We now distinguish three cases.
 
 \begin{enumerate}
\item Suppose first $m = n$. By the above, $\sg(K_n \cp K_n) \geq n$. Take $n$ vertices, one in each row and one in each column. Since $\diam(K_n \cp K_n) = 2$, these $n$ vertices can cover at most $\binom{n}{2}$ other vertices of $K_n \cp K_n$. Moreover, at most one row and at most one column can be covered completely with geodesics between them. Hence, at least $\binom{n}{2}$ vertices of $K_n \cp K_n$ remain uncovered. As at least $n-1$ rows and columns are still uncovered, it follows that at least $n-1$ more vertices are needed to cover them. Therefore, $\sg(K_n \cp K_n) \geq n + (n-1) = 2n - 1$. 

Consider the set $S = S_1 \cup S_2$ where $S_1 =  \{ (i, i) \; ; \; i \in [n] \}$ and $S_2 = \{ (i, i+1) \; ; \; i \in [n-1] \}$ (cf.~Fig.~\ref{fig:complete1}). 

% \iffalse
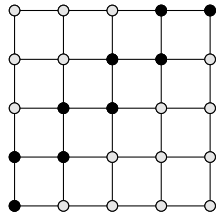
\begin{figure}[!ht]
\begin{center}
    \begin{tikzpicture}[scale=0.65]
    \pgfmathtruncatemacro{\m}{5}
    \pgfmathtruncatemacro{\n}{5}
    \pgfmathtruncatemacro{\nm}{4}
    \begin{scope}
        \foreach \x in {1,...,\m}
            \foreach \y in {1,...,\n}
                \node (\x,\y) at (\x, \y) {};
        \foreach \x [remember=\x as \lastx (initially 1)] in {1,...,\m}
            \foreach \y in {1,...,\n}
                \path (\x, \y) edge (\lastx, \y);
                
        \foreach \y [remember=\y as \lasty (initially 1)] in {1,...,\n}
            \foreach \x in {1,...,\m}
                \path (\x, \y) edge (\x, \lasty);
        \foreach \x in {1,...,\m}
            \foreach \y in {1,...,\n}
                \node (\x,\y) at (\x, \y) {};
                
        \foreach \x in {1,...,\n}
            \node[fill=black] (\x) at (\x, \x) {};
        \foreach \x in {1,...,\nm}
            \node[fill=black] (\x) at (\x, \x+1) {}; 
    \end{scope}
    \end{tikzpicture}
    \caption{A strong geodetic set of $K_5 \cp K_5$.}
    \label{fig:complete1}
\end{center}
\end{figure}
% \fi

Fix geodesics for $\widetilde{I}(S)$ in such a way that geodesics between vertices from $S_1$ cover all the vertices $\{ (i, j) \; ; \; i \geq j \}$ and geodesics between vertices from $S_2$ cover the vertices $\{ (i, j) \; ; \; i < j \}$. Thus $S$ is a strong geodetic set of size $2n -1$. Hence, $\sg(K_n \cp K_n) = 2n-1$.  

\item Suppose next $n < m < 2n$. Consider an arbitrary strong geodetic set $S'$ of $K_m\cp K_n$. Since $S'$ contains at least one vertex from each row and at least one vertex from each column, we may without loss of generality assume that $S_1\cup S_2\subseteq S'$, where $S_1 = \{ (i, i) \; ; \; i \in [n] \}$ and  $S_2 = \{ (n + i, i) \; ; \; i \in [m-n] \}$. Consider the disjoint sets 
\begin{eqnarray*}
A & = & [m-n] \times \{ m-n+1, \ldots, n \}, \\
B & = & \{ n+1, \ldots, m \} \times \{ m-n+1, \ldots, n \}, \\
C & = & \{ m-n+1, \ldots, n \} \times [m-n],\\
D & = & \{ (i, j) \; ; \; i < j, i, j \in \{ m-n+1, \ldots, n \} \},\quad \mbox{and} \\
E & = & \{ (i, j) \; ; \; i > j, i, j \in \{ m-n+1, \ldots, n \} \},
\end{eqnarray*}
which are shown in Fig.~\ref{fig:complete_sets} for the case $K_{10} \cp K_7$.

\begin{figure}[!ht]
\begin{center}
    \begin{tikzpicture}[scale=0.65]
    \pgfmathtruncatemacro{\m}{10}
    \pgfmathtruncatemacro{\n}{7}
    \pgfmathtruncatemacro{\mn}{3}
    \pgfmathtruncatemacro{\mnp}{4}
    \begin{scope}
        \foreach \x in {1,...,\m}
            \foreach \y in {1,...,\n}
                \node (\x,\y) at (\x, \y) {};
        \foreach \x [remember=\x as \lastx (initially 1)] in {1,...,\m}
            \foreach \y in {1,...,\n}
                \path[black!70] (\x, \y) edge (\lastx, \y);
                
        \foreach \y [remember=\y as \lasty (initially 1)] in {1,...,\n}
            \foreach \x in {1,...,\m}
                \path[black!70] (\x, \y) edge (\x, \lasty);
        \foreach \x in {1,...,\m}
            \foreach \y in {1,...,\n}
                \node (\x,\y) at (\x, \y) {};
                
        \foreach \x in {1,...,\n}
            \node[fill=black!70] (\x) at (\x, \x) {};

        \foreach \x in {1,...,\mn}
            \node[fill=black!70] (\x') at (\n + \x, \x) {};
            
        \path[draw, black] (1 - 0.2, \mnp - 0.2) -- (\mn + 0.2, \mnp - 0.2) -- (\mn + 0.2, \n + 0.2) -- (1 - 0.2, \n + 0.2) -- (1 - 0.2, \mnp - 0.2);
        \node[draw=none, fill=none] at (1.5, \n-0.5) {$A$};
        \path[draw, black] (\n + 1 - 0.2, \mnp - 0.2) -- (\m + 0.2, \mnp - 0.2) -- (\m + 0.2, \n + 0.2) -- (\n + 1 - 0.2, \n + 0.2) -- (\n + 1 - 0.2, \mnp - 0.2);
        \node[draw=none, fill=none] at (\m - 0.5, \n-0.5) {$B$};
        \path[draw, black] (\mnp - 0.2, 1- 0.2) -- (\n + 0.2, 1 - 0.2) -- (\n + 0.2, \mn + 0.2) -- (\mnp - 0.2, \mn + 0.2) -- (\mnp - 0.2, 1 - 0.2);
        \node[draw=none, fill=none] at (\mnp + 1 + 0.5, 1 + 0.5) {$C$};
        \path[draw, black] (\mnp - 0.2, \mnp+1 - 0.4) -- (\n - 1 + 0.4, \n + 0.2) -- (\mnp - 0.2, \n + 0.2) -- (\mnp - 0.2, \mnp+1 - 0.4);
        \node[draw=none, fill=none] at (\mnp + 0.5, \n-0.5) {$D$};
        \path[draw, black] (\mnp+1-0.4, \mnp-0.2) -- (\n+0.2, \mnp-0.2) -- (\n+0.2, \n-1 + 0.4) -- (\mnp+1-0.4, \mnp-0.2);
        \node[draw=none, fill=none] at (\n - 0.5, \mnp + 0.5) {$E$};
    \end{scope}
    \end{tikzpicture}
    \caption{Sets $A, B, C, D$ and $E$ of $K_{10} \cp K_7$.}
    \label{fig:complete_sets}
\end{center}
\end{figure}
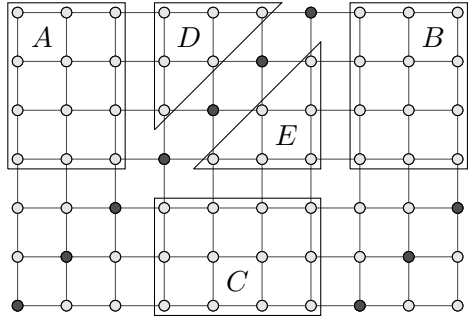

Vertices in $A$ can only be covered with geodesics between vertices from $S_1$, thus these geodesics cannot cover $C$. The set $B$ can only be covered with geodesics between vertices from $S_1$ and $S_2$ and thus these geodesics cannot cover $C$. Hence, $C$ is left uncovered. Similarly we observe, that either $D$ or $E$ is left uncovered. It follows that vertices lying in $2n - m$ different columns and vertices from $n-1$ different rows are left uncovered. To cover them, at least $\min\{2n - m, n-1\}$ additional vertices must be added to $S'$. As $m > n$, we have $\min\{2n - m, n-1\} = 2n - m$. Hence, $\sg(K_m \cp K_n) \geq m + (2n - m) = 2n$.

Consider the set $S = S_1 \cup S_2 \cup S_3$, where $S_1$ and $S_2$ are as above and $S_3 = \{ (i, 1) \; ; \; i \in \{ m-n+1, \ldots, n \} \}$ (cf.~Fig.~\ref{fig:complete2}). Denote $S_1 = S_1^d \cup S_1^u$, where $S_1^d = \{ (i, i) \; ; \; i \in [m-n] \}$ and $S_1^u = \{ (i, i) \; ; \; i \in \{ m-n+1, \ldots, n \} \}$.  

% \iffalse
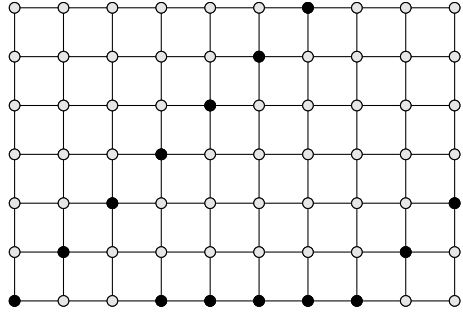
\begin{figure}[!ht]
\begin{center}
    \begin{tikzpicture}[scale=0.65]
    \pgfmathtruncatemacro{\m}{10}
    \pgfmathtruncatemacro{\n}{7}
    \pgfmathtruncatemacro{\mn}{3}
    \pgfmathtruncatemacro{\mnp}{4}
    \begin{scope}
        \foreach \x in {1,...,\m}
            \foreach \y in {1,...,\n}
                \node (\x,\y) at (\x, \y) {};
        \foreach \x [remember=\x as \lastx (initially 1)] in {1,...,\m}
            \foreach \y in {1,...,\n}
                \path (\x, \y) edge (\lastx, \y);
                
        \foreach \y [remember=\y as \lasty (initially 1)] in {1,...,\n}
            \foreach \x in {1,...,\m}
                \path (\x, \y) edge (\x, \lasty);
        \foreach \x in {1,...,\m}
            \foreach \y in {1,...,\n}
                \node (\x,\y) at (\x, \y) {};
                
        \foreach \x in {1,...,\n}
            \node[fill=black] (\x) at (\x, \x) {};
            
        \foreach \x in {\mnp,...,\n}
            \node[fill=black] (\x) at (\x, 1) {};
            
        \foreach \x in {1,...,\mn}
            \node[fill=black] (\x') at (\n + \x, \x) {};            
    \end{scope}
    \end{tikzpicture}
    \caption{A strong geodetic set of $K_{10} \cp K_7$.}
    \label{fig:complete2}
\end{center}
\end{figure}
% \fi

Fix geodesics between vertices in $S_1$ to cover $\{ (i, j) \; ; \; i < j, i, j \in [n] \}$, geodesics between vertices in $S_2$ to cover $\{ (i, j) \; ; \; i < j, i \in \{ n+1, \ldots, m \}, j \in [m-n] \}$, geodesics between $S_1^d$ and $S_2$ to cover $\{ (i, j) \; ; \; i > j, i \in [m-n] \cup \{n+1, \ldots, m\}, j \in [m-n] \}$ and geodesics between $S_1^u$ and $S_2$ to cover $\{ n+1, \ldots, m \} \times \{ m-n+1, \ldots, n \}$. Additionaly, fix geodesics $(v, 1) \sim (v, i) \sim (i, i)$ for each $v \in S_3$ and $i \in [n]$. Now it is clear that $S$ is a strong geodetic set of size $2n$. Hence, $\sg(K_m \cp K_n) = 2n$. 

\item Suppose finally $m \geq 2n$. We already know that $\sg(K_m \cp K_n) \geq m$. Define $S = S_l \cup S_m \cup S_r$, where 
\begin{eqnarray*}
S_l & = & \{ (i, i) \; ; \; i \in [n] \}, \\
S_m & = & \{ (i, 1) \; ; \; i \in \{ n+1, \ldots, m-n \} \}, \quad \mbox{and} \\
S_r & = & \{ (m-n + i, i) \; ; \; i \in [n] \}, 
\end{eqnarray*}
cf.~Fig.~\ref{fig:complete3}, where $S$ is shown for the case $K_{12}\cp K_4$.

% \iffalse
\begin{figure}[!ht]
\begin{center}
    \begin{tikzpicture}[scale=0.65]
    \pgfmathtruncatemacro{\m}{12}
    \pgfmathtruncatemacro{\n}{4}
    \pgfmathtruncatemacro{\np}{5}
    \pgfmathtruncatemacro{\mn}{8}
    \begin{scope}
        \foreach \x in {1,...,\m}
            \foreach \y in {1,...,\n}
                \node (\x,\y) at (\x, \y) {};
        \foreach \x [remember=\x as \lastx (initially 1)] in {1,...,\m}
            \foreach \y in {1,...,\n}
                \path (\x, \y) edge (\lastx, \y);
                
        \foreach \y [remember=\y as \lasty (initially 1)] in {1,...,\n}
            \foreach \x in {1,...,\m}
                \path (\x, \y) edge (\x, \lasty);
        \foreach \x in {1,...,\m}
            \foreach \y in {1,...,\n}
                \node (\x,\y) at (\x, \y) {};
                
        \foreach \x in {1,...,\n}
            \node[fill=black] (\x) at (\x, \x) {};
            
        \foreach \x in {\np,...,\mn}
            \node[fill=black] (\x) at (\x, 1) {};
            
        \foreach \x in {1,...,\n}
            \node[fill=black] (\x') at (\mn + \x, \x) {};            
    \end{scope}
    \end{tikzpicture}
    \caption{A strong geodetic set of $K_{12} \cp K_4$.}
    \label{fig:complete3}
\end{center}
\end{figure}
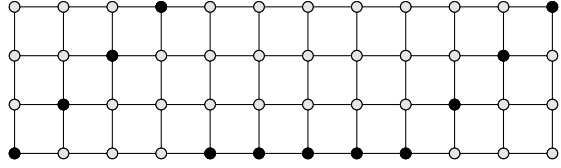
% \fi

Fix geodesics between vertices from $S_l$ to cover vertices $\{ (i, j) \; ; \; i \geq j, i,j \in [n] \}$, geodesics between vertices from $S_r$ to cover $\{ (m - n + i, j) \; ; \; i \geq j, i,j \in [n] \}$, geodesics between sets $S_l$ and $S_r$ to cover $\{ (i, j) \; ; \; i \leq j, i,j \in [n] \} \cup \{ (m - n + i, j) \; ; \; i \leq j, i,j \in [n] \}$ and geodesics between a vertex $v \in S_m$ and vertices from $S_l$ to cover $\{ (v, i) \; ; \; i \in [n] \}$. Hence $S$ is a strong geodetic set of $K_m \cp K_n$ and $|S| = m$. We conclude that $\sg(K_m \cp K_n) = m$. \qed
\end{enumerate}

From Theorem~\ref{thm:KmKn} we infer that among Cartesian products of complete graphs the upper bound of Theorem~\ref{thm:generalUpper} is sharp only for $K_1 \cp K_1$, $K_2 \cp K_2$, and $K_3 \cp K_2$.

\medskip
Until now we have considered general upper bounds on $\sg(G\cp H)$ and obtained several exact values. Hence it would also be of interest to have some general lower bound(s). For this sake we pose: 

\begin{conjecture}
\label{conj:lower-bound}
If $G$ is a graph with $\n(G)\ge 2$, then $\sg(G \cp K_2) \geq \sg(G)$.
\end{conjecture}

\begin{figure}[!ht]
\begin{center}
    \begin{tikzpicture}
    \pgfmathtruncatemacro{\n}{5}
    \pgfmathtruncatemacro{\r}{1}
    \begin{scope}
        \foreach \x in {0,...,\n}
            \node (\x) at (\x*360/\n:\r cm) {};
        \foreach \x in {0,...,3}
            \path (\x) edge (4);
    \end{scope}
    
    \begin{scope}[xshift = 2.5cm]
        \foreach \x in {0,...,\n}
            \node (\x) at (\x*360/\n:\r cm) {};
        \foreach \x in {0,...,3}
            \path (\x) edge (4);
        \path (0) edge (3);
    \end{scope}
    
    \begin{scope}[xshift = 5cm]
        \foreach \x in {0,...,\n}
            \node (\x) at (\x*360/\n:\r cm) {};
        \foreach \x in {0,...,3}
            \path (\x) edge (4);
        \path (0) edge (2);
        \path (1) edge (3);
    \end{scope}
    
    \begin{scope}[xshift = 7.5cm]
        \foreach \x in {0,...,\n}
            \node (\x) at (\x*360/\n:\r cm) {};
        \foreach \x in {0,...,3}
            \path (\x) edge (4);
        \path (0) edge (2);
        \path (0) edge (3);
        \path (2) edge (3);
    \end{scope}
    
    \begin{scope}[xshift = 1.25cm, yshift = -2.5cm]
        \foreach \x in {0,...,\n}
            \node (\x) at (\x*360/\n:\r cm) {};
        \foreach \x in {0,...,3}
            \path (\x) edge (4);
        \path (0) edge (2);
        \path (0) edge (3);
        \path (2) edge (1);
        \path (3) edge (1);
    \end{scope}
    
    \begin{scope}[xshift = 3.75cm, yshift = -2.5cm]
        \foreach \x in {0,...,\n}
            \node (\x) at (\x*360/\n:\r cm) {};
        \foreach \x in {0,...,3}
            \path (\x) edge (4);
        \path (0) edge (2);
        \path (0) edge (3);
        \path (2) edge (1);
        \path (3) edge (1);
        \path (3) edge (2);
    \end{scope}
    
    \begin{scope}[xshift = 6.25cm, yshift = -2.5cm]
        \foreach \x in {0,...,\n}
            \node (\x) at (\x*360/\n:\r cm) {};
        \foreach \x in {0,...,3}
            \path (\x) edge (4);
        \path (0) edge (2);
        \path (0) edge (3);
        \path (2) edge (1);
        \path (3) edge (1);
        \path (3) edge (2);
        \path (0) edge (1);
    \end{scope}
    
    \end{tikzpicture}
    \caption{Graphs $G$ on five vertices with $\sg(G) = \sg(G \cp K_2)$.}
    \label{fig:n=5}
\end{center}
\end{figure}
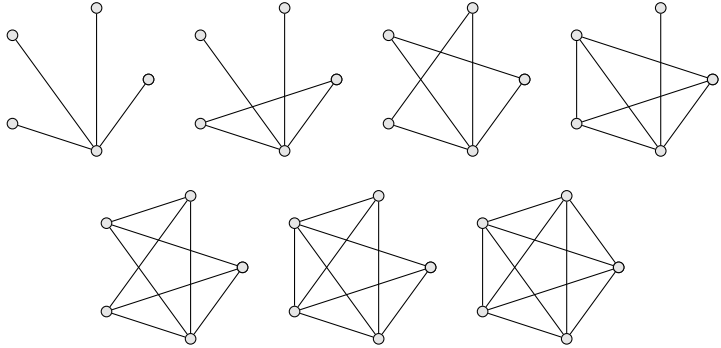

If Conjecture~\ref{conj:lower-bound} is true, then it is best possible as demonstrated by Theorem~\ref{thm:equality_prisms}. We have also verified the cojecture by computer for all graphs $G$ with $\n(G) \leq 7$. The equality is never attained for $\n(G) \le 3$. For $\n(G) = 4$ the only equality case is $G = K_4$, while for $\n(G) = 5$ and $6$ there are more equality cases. For $\n(G) = 5$ all of them are shown in Fig.~\ref{fig:n=5}. For $\n(G)=6$ the variety of equality graphs is too large to be drawn here.

More generally as Conjecture~\ref{conj:lower-bound}, we pose the following

\begin{problem}
\label{problem:lower-bound}
Is it true that if $G$ and $H$ are graphs, then $\sg(G \cp H) \geq \max\{ \sg(G), \sg(H) \}$?
\end{problem}

Again, if the answer to Problem~\ref{problem:lower-bound} is positive, then the result is best possible as demonstrated by Proposition~\ref{prop:starPn} and by Theorem~\ref{thm:KmKn} for $m\ge 2n$. 

%%%%%%%%%%%%%%%%%%%%%%%%%%%%%%%%%%%%%%%%%%%%%%%%%%%%%%%%%%
%%%%%%%%%%%%%%%%%%%%%%%%%%%%%%%%%%%%%%%%%%%%%%%%%%%%%%%%%%
\section{The strong geodetic number of subgraphs}
\label{sec:subgraphs}
%%%%%%%%%%%%%%%%%%%%%%%%%%%%%%%%%%%%%%%%%%%%%%%%%%%%%%%%%%
%%%%%%%%%%%%%%%%%%%%%%%%%%%%%%%%%%%%%%%%%%%%%%%%%%%%%%%%%%

Since layers of Cartesian products are subgraphs that possess several distinguishing properties, a way to attack Conjecture~\ref{conj:lower-bound} would be to understand the relation between the strong geodetic number of a graph and its subgraphs. This is a fundamental question for any graph invariant and has not yet been studied for the strong geodetic number. The main message of this section is that in general there is no such relation, even for subgraphs with a very special structure such as layers in products.
  
\subsection*{Induced subgraphs}

First we observe that there is no connection between a strong geodetic number of a graph and a strong geodetic number of its (induced) subgraph. 

\begin{figure}[!ht]
\begin{center}
    \begin{tikzpicture}
    \pgfmathtruncatemacro{\n}{8}
    \begin{scope}
        \foreach \x in {1,...,\n}
            \node (\x) at (\x, 0) {};
        \foreach \x [remember=\x as \lastx (initially 1)] in {1,...,\n}
            \path (\x) edge (\lastx);
            
        \foreach \x in {1,...,\n}
            \node (\x') at (\x, 1) {};
        \foreach \x [remember=\x as \lastx (initially 1)] in {1,...,\n}
            \path (\x') edge (\lastx');
            
        \foreach \x in {1,...,\n}
            \path (\x) edge (\x');
        
        \node[fill=black] (1) at (1,0) {};
        \node[fill=black] (\n) at (\n,0) {};
        \node[fill=black] (\n') at (\n,1) {};
    \end{scope}
    
    \begin{scope}[xshift=7.5cm]
        \foreach \x in {1,3,5,7}
            \node[fill=black] (\x) at (\x, 0) {};
            
        \foreach \x in {1,...,\n}
            \node (\x') at (\x, 1) {};
        \foreach \x [remember=\x as \lastx (initially 1)] in {1,...,\n}
            \path (\x') edge (\lastx');
            
        \foreach \x in {1,3,5,7}
            \path (\x) edge (\x');
        
        \node[fill=black] (\n') at (\n,1) {};
    \end{scope}
    \end{tikzpicture}
    \caption{The strong geodetic sets of graphs $G_4$ and its subgraph $H_4$.}
    \label{fig:subgraph}
\end{center}
\end{figure}
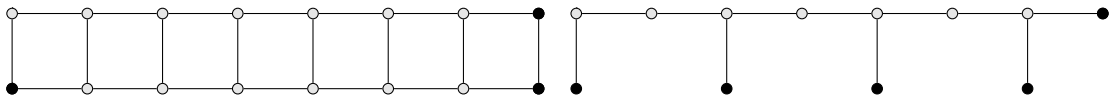

Let $G_n = P_{2n} \cp K_2$ and $H_n$ its subgraph induced on vertices $V(G_n) - \{ (2i, 1) \; ; \; i \in [n] \}$ (cf.~Fig.~\ref{fig:subgraph}). Clearly, $\sg(G_n) = 3$, as $\{(1, 1), (2n, 1), (2n, 2)\}$ is a strong geodetic set. The subgraph $H_n$ is a tree with $n+1$ leaves, thus $\sg(H_n) = n+1$. Hence, the strong geodetic number of an induced subgraph can be arbitrarily larger that the strong geodetic number of a graph. The converse is also true. Consider $H = P_n$ as a(n) (induced) subgraph of some tree $T$. It holds $\sg(H) = 2$, but the strong geodetic number of $T$ can be arbitrarily large (and equals the number of its leaves). 

% \iffalse
% \fi

\subsection*{Convex subgraphs}

A subgraph $H$ of graph $G$ is \emph{convex} if every shortest path in $G$ between vertices from $H$ lies entirely in $H$. This is a stronger concept than induced subgraphs. Layers of Cartesian products are convex. 

As paths are convex subgraphs of trees, it is clear that the strong geodetic number of a graph can be arbitrarily larger than the strong geodetic number of its convex subgraphs. The following example shows that the converse also holds. 

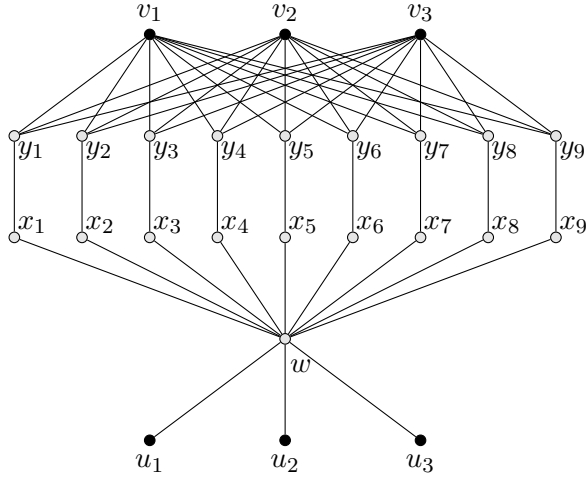
\begin{figure}[!ht]
\begin{center}
    \begin{tikzpicture}[scale=0.9]
    \pgfmathtruncatemacro{\k}{3}
    \pgfmathtruncatemacro{\l}{3}
    \pgfmathtruncatemacro{\kl}{9}

    \begin{scope}
        \foreach \x in {1,...,\k}
            \node[fill=black, label=below: $u_\x$] (u_\x) at (1+2*\x, -1) {};
            
        \node[label={[label distance=5pt]290:$w$}] (w) at (5, 0.5) {};
            
        \foreach \x in {1,...,\kl}
            \node[label=5:$x_\x$] (x_\x) at (\x, 2) {};
            
        \foreach \x in {1,...,\kl}
            \node[label=330:$y_\x$] (y_\x) at (\x, 3.5) {};
            
        \foreach \x in {1,...,\l}
            \node[fill=black, label=above:$v_\x$] (v_\x) at (1+2*\x, 5) {};
                    
        \foreach \x in {1,...,\k}
            \path (u_\x) edge (w);
            
        \foreach \x in {1,...,\kl}
            \path (x_\x) edge (w);
            
        \foreach \x in {1,...,\kl}
            \path (x_\x) edge (y_\x);
            
        \foreach \x in {1,...,\l}
            \foreach \y in {1,...,\kl}
                \path (v_\x) edge (y_\y);
    \end{scope}
    \end{tikzpicture}
    \caption{The graph $G^c_{3,3}$.}
    \label{fig:convex}
\end{center}
\end{figure}

Let $k, l \in \N$. Define $G^c_{k, l}$ to be the graph with $V(G^c_{k,l}) = \{ u_1, \ldots, u_k \} \cup \{w\} \cup \{x_1, y_1, \ldots, x_{kl}, y_{kl} \} \cup \{ v_1, \ldots, v_l \}$ and edges $w \sim u_i$ for $i \in [k]$, $w \sim x_i$ for $i \in [kl]$, $x_i \sim y_i$ for $i \in [kl]$ and $y_i \sim v_j$ for all $i \in [kl]$ and $j \in [l]$ (cf.~Fig.~\ref{fig:convex}). Let $H$ be its subgraph induced by $\{w\} \cup \{ x_1, \ldots, x_{kl} \}$. Note that $H$ is a convex subgraph with $\sg(H) = kl$ (as it is a tree).

As vertices $\{u_1, \ldots, u_k\}$ are simplicial, they lie in any strong geodetic set of $G^c_{k,l}$. But due to the structure of the graph, each vertex $v_i$ must also lie in any strong geodetic set. Hence, $\sg(G^c_{k, l}) \geq k + l$. Consider the set $S = \{u_1, \ldots, u_k\} \cup \{ v_1, \ldots, v_l \}$ and fix the geodesics $u_i \sim w \sim x_{(i-1) l + j} \sim y_{(i-1) l + j} \sim v_j$ for all $i \in [k]$ and $j \in [l]$. These geodesics cover all vertices of a graph, hence $\sg(G^c_{k,l}) = k+l$, which is arbitrarily smaller that $kl$, the strong geodetic number of the convex subgraph $H$.

\subsection*{Gated subgraphs}

A subgraph $H$ of graph $G$ is \emph{gated} if for every $v \in V(G)$ there exists an $x \in V(H)$ that lies on a shortest $u, v$-path for every $u \in V(H)$. Every gated subgraph is convex~\cite{handbook}. Layers of Cartesian product are not only convex but also gated. 

Unfortunately, there is also no connection between the strong geodetic number of a graph and its gated subgraphs. Again, as paths are gated subgraphs of trees, the strong geodetic number of a graph can be arbitrarily larger than the strong geodetic number of its gated subgraphs. The following example shows that the converse is also true.

\begin{figure}[!ht]
\begin{center}
    \begin{tikzpicture}[scale=0.9]
    \pgfmathtruncatemacro{\k}{2}
    \pgfmathtruncatemacro{\l}{3}
    \pgfmathtruncatemacro{\kl}{6}

    \begin{scope}
        \foreach \x in {1,...,\k}
            \node[fill=black, label=left:$x_\x$] (x_\x) at (0, \x+2) {};
        \node[label={[label distance=3pt]270:$x$}] (x) at (1.5, \kl/2 + 1/2) {};
        \foreach \x in {1,...,\k}
            \foreach \y in {1,...,\l}
                \node[label=below:$v_{\x, \y}$] (v_\x_\y) at (3, \y + 3*\x - 3) {};
        \node[label={[label distance=3pt]270:$y$}] (y) at (4.5, \kl/2 + 1/2) {};
        \foreach \x in {1,...,\l}
            \node[fill=black, label=right:$y_\x$] (y_\x) at (6, \x+3/2) {};
            
        \foreach \x in {1,...,\k}
            \path (x_\x) edge (x);
        \foreach \x in {1,...,\l}
            \path (y_\x) edge (y);
        \foreach \x in {1,...,\k}
            \foreach \y in {1,...,\l}
                \path (x) edge (v_\x_\y);
        \foreach \x in {1,...,\k}
            \foreach \y in {1,...,\l}
                \path (y) edge (v_\x_\y);
        
    \end{scope}
    \end{tikzpicture}
    \caption{The graph $G^g_{2,3}$.}
    \label{fig:gated}
\end{center}
\end{figure}
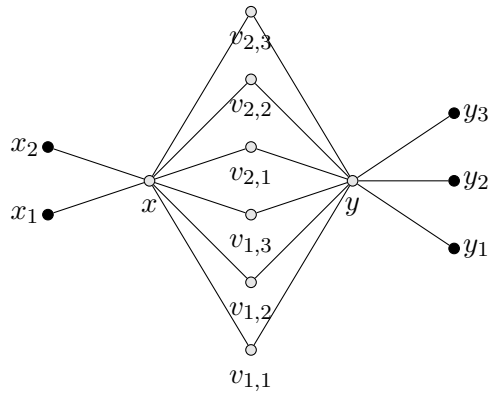

Let $k, l \in \N$ such that $kl \geq 5$. Define the graph $G^g_{k, l}$ with vertices $\{ x, y \} \cup \{ v_{i, j} \; ; \; i \in [k], j \in [l] \} \cup \{ x_1, \ldots, x_k \} \cup \{ y_1, \ldots, y_l \}$ and edges $x \sim x_i$ for $i \in [k]$, $y \sim y_j$ for $j \in [l]$, $x \sim v_{i, j} \sim y$ for $i \in [k], j\in [l]$ (cf.~Fig.~\ref{fig:gated}).   

% \iffalse
% \fi

Let $S = \{ x_1, \ldots, x_k, y_1, \ldots, y_l \}$. Vertices in $S$ are all simplicial, thus $\sg(G^g_{k,l}) \geq |S| = k + l$. If we fix geodesics $x_i \sim x \sim v_{i,j} \sim y \sim y_j$ for all $i \in [k], j \in [l]$, then it is clear that $S$ is a strong geodetic set. Hence, $\sg(G^g_{k,l}) = k + l$.

Let $H$ be a subgraph of $G$ induced on the vertex set $\{ x, y \} \cup \{ v_{i, j} \; ; \; i \in [k], j \in [l] \}$. Clearly, $H \cong K_{2,kl}$. A subgraph $H$ is gated in $G$. It follows from $kl \geq 5$, that $\binom {kl - 1}{2} \geq kl$ and thus by~\cite{bipartite} it holds that $\sg(H) = kl$. Hence, the strong geodetic number of a gated subgraph can be arbitrarily larger than the strong geodetic number of a graph.

\end{document}